\documentclass[italian,UKenglish,12pt]{amsart}
\usepackage{etex}
\usepackage[utf8]{inputenc}
\usepackage[a4paper,lmargin=2.8cm,rmargin=2.8cm,tmargin=3.0cm,bmargin=3.0cm]{geometry}
\usepackage{multirow}

\usepackage{babel}

\usepackage{rotating,color,hyperref,appendix}
\usepackage{tocvsec2}
\usepackage{booktabs}
\usepackage{float}
\usepackage{tikz}
\usetikzlibrary{positioning,arrows,shapes,decorations.pathmorphing,shadows}

\usepackage{amsfonts}

\newcommand{\Spin}[1]{\ensuremath{\text{\upshape\rmfamily Spin}(#1)}}

\newcommand{\Sp}[1]{\mathrm{Sp}(#1)}
\newcommand{\SO}[1]{\mathrm{SO}(#1)}
\newcommand{\SU}[1]{\mathrm{SU}(#1)}
\newcommand{\U}[1]{\mathrm{U}(#1)}

\newcommand{\spinform}[1]{\ensuremath{\Phi_{\Spin{#1}}}}

\newcommand{\LH}{L^\HH}

\newcommand{\liespin}[1]{\mathop{\mathfrak{spin}}(#1)}
\newcommand{\liesu}[1]{\mathop{\mathfrak{su}}(#1)}

\newcommand{\liesp}[1]{\mathop{\mathfrak{sp}}(#1)}
\newcommand{\lieso}[1]{\mathop{\mathfrak{so}}(#1)}

\newcommand{\Id}{\text{Id}}

\newcommand{\End}[1]{\mathrm{End}(#1)}

\newcommand{\CC}{\mathbb{C}}
\newcommand{\HH}{\mathbb{H}}
\newcommand{\RR}{\mathbb{R}}
\newcommand{\ZZ}{\mathbb{Z}}

\newcommand{\OO}{\mathbb{O}}

\newcommand{\Gtwo}{\mathrm{G}_2}
\newcommand{\EI}{\mathrm{E\,I}}
\newcommand{\EII}{\mathrm{E\,II}}
\newcommand{\EIII}{\mathrm{E\,III}}
\newcommand{\EIV}{\mathrm{E\,IV}}
\newcommand{\EV}{\mathrm{E\,V}}
\newcommand{\EVI}{\mathrm{E\,VI}}
\newcommand{\EVII}{\mathrm{E\,VII}}
\newcommand{\EVIII}{\mathrm{E\,VIII}}
\newcommand{\EIX}{\mathrm{E\,IX}}
\newcommand{\FI}{\mathrm{F\,I}}
\newcommand{\FII}{\mathrm{F\,II}}
\newcommand{\GI}{\mathrm{G\,I}}
\newcommand{\Esix}{\mathrm{E}_6}
\newcommand{\Eseven}{\mathrm{E}_7}
\newcommand{\Eeight}{\mathrm{E}_8}
\newcommand{\Ffour}{\mathrm{F}_4}

\numberwithin{equation}{section}

\theoremstyle{plain}
\newtheorem{te}{Theorem}[section]
\newtheorem*{te*}{Theorem}
\newtheorem{pr}[te]{Proposition}

\theoremstyle{definition}
\newtheorem{de}[te]{Definition}

\theoremstyle{remark}
\newtheorem{re}[te]{Remark}

\date{\today}


\begin{document}

\begin{otherlanguage}{italian}
\author{Paolo Piccinni}
\address{Sapienza-Universit\`a di Roma\\ Dipartimento di Matematica, piazzale Aldo Moro 2, I-00185, Roma, Italy}
\email{piccinni@mat.uniroma1.it}
\end{otherlanguage}

\title[On the cohomology]{On the cohomology\\ of some exceptional symmetric spaces}

\thanks{The author was supported by the GNSAGA group of INdAM and by the research project "Polynomial identities and combinatorial methods in algebraic and geometric structures" of Sapienza Universit\`a di Roma}

\vspace{1.5cm}
\maketitle

%
\begin{center}
\emph{Dedicated to Simon Salamon on the occasion of his 60th birthday}
\end{center}


\begin{abstract}
This is a survey on the construction of a canonical or "octonionic K\"ahler" 8-form, representing one of the generators of the cohomology of the four Cayley-Rosenfeld projective planes. The construction, in terms of the associated even Clifford structures, draws a parallel with that of the quaternion K\"ahler 4-form. We point out how these notions allow to describe the primitive Betti numbers with respect to different even Clifford structures, on most of the exceptional symmetric spaces of compact type. \end{abstract}

\vspace{1cm}

\section{Introduction}
\noindent The exceptional Riemannian symmetric spaces of compact type 
\[
\EI, \, \EII, \; \EIII, \; \EIV, \; \EV, \; \EVI, \; \EVII, \; \EVIII, \; \EIX, \; \FI, \; \FII,  \; \GI
\] 
are part of the E. Cartan classification. 

Among them, the two Hermitian symmetric spaces 
\[
\EIII = \frac{\Esix}{\mathrm{Spin}(10)\cdot\mathrm{U}(1)} \quad \text{and} \quad \EVII =  \frac{\Eseven}{\Esix \cdot\mathrm{U}(1)}
\]
are certainly notable. As Fano manifolds, they can be realized as smooth complex projective varieties. As such, $\EIII$ is also called the \emph{fourth Severi variety}, a complex 16-dimensional projective variety in $\CC P^{26}$, characterized as one of the four smooth projective varieties of small critical codimension in their ambient $\CC P^{N}$, and that are unable to fill it through their secant and tangent lines \cite{ZakSeV}. The projective model of $\EVII$ is instead known as the \emph{Freudenthal variety}, a complex 27-dimensional projective variety in $\CC P^{55}$, considered for example in the sequel of papers \cite{cmp}. 

Next, among the listed symmetric spaces, the five Wolf spaces
\[
\begin{split}
\EII= \frac{\Esix}{\mathrm{SU}(6)\cdot\mathrm{Sp}(1)},\quad  \EVI=  \frac{\Eseven}{\mathrm{Spin}(12)\cdot\mathrm{Sp}(1)},  \quad & \EIX = \frac{\Eeight}{\Eseven\cdot\mathrm{Sp}(1)}, \\ 
\FI =  \frac{\Ffour}{\mathrm{Sp}(3)\cdot\mathrm{Sp}(1)},  \quad  \GI = \frac{\Gtwo}{\mathrm{SO}(4)}&
\end{split}
\]
give evidence for the long lasting LeBrun-Salamon conjecture, being the only known sporadic examples of positive quaternion K\"ahler manifolds.

Thus, seven of the twelve exceptional Riemannian symmetric spaces of compact type are either K\"ahler or quaternion K\"ahler. Accordingly, one of their de Rham cohomology generators is represented by a K\"ahler 2-form or a quaternion K\"ahler 4-form, and any further cohomology generators can be looked as primitive in the sense of the Lefschetz decomposition.

The notion of \emph{even Clifford structure}, introduced some years ago by A. Moroianu and U. Semmelmann \cite{ms}, allows not only to deal simultaneously with K\"ahler and quaternion K\"ahler manifolds, but also to recognize further interesting geometries fitting into the notion. Among them, and just looking at the exceptional Riemannian symmetric spaces of compact type, there are even Clifford structures, related with octonions, on the following \emph{Cayley-Rosenfeld projective planes}
\[
\begin{split}
\EIII = \frac{\Esix}{\mathrm{Spin}(10)\cdot\mathrm{U}(1)}, \qquad & \EVI=  \frac{\Eseven}{\mathrm{Spin}(12)\cdot\mathrm{Sp}(1)}, \\ \EVIII = \frac{\Eeight}{\mathrm{Spin}^+(16)}, \qquad & \FII =  \frac{\Ffour}{\mathrm{Spin}(9)}.
\end{split}
\]

An even Clifford structure is defined as the datum, on a Riemannian manifold $(M,g)$, of a real oriented Euclidean vector bundle $(E,h)$, together with an algebra bundle morphism 
\[
\varphi: \; \text{Cl}^0(E) \rightarrow \End{TM}
\] 
mapping $\Lambda^2 E$ into skew-symmetric endomorphisms. The rank $r$ of $E$ is said to be the \emph{rank of the even Clifford structure}. One easily recognizes that K\"ahler and quaternion K\"ahler metrics correspond to a choice of such a vector bundle $E$ with $r=2,3$ respectively, and that for the four Cayley-Rosenfeld projective planes $\EIII, \, \EVI, \, \EVIII , \, \FII$ there is a similar vector bundle $E$ with $r=10,12,16,9$, cf. \cite{ms}. 

Thus, among the exceptional symmetric spaces of compact type, there are two spaces admitting two distinct even Clifford structures. Namely, the Hermitian symmetric $\EIII$ has both a rank $2$ and a rank $10$ even Clifford structure, and the quaternion K\"ahler $\EVI$ has both a rank $3$ and a rank $12$ even Clifford structure. Moreover, all the even Clifford structures we are here considering on our symmetric spaces are \emph{parallel}, i. e. there is a metric connection $\nabla^E$ on $(E,h)$ such that $\varphi$ is connection preserving: $$\varphi(\nabla^E_X \sigma)=\nabla^g_X\varphi(\sigma),$$ for every tangent vector $X \in TM$, section $\sigma$ of $\text{Cl}^0(E)$, and where $\nabla^g$ is the Levi Civita connection of the Riemannian metric $g$. For simplicity, we will call \emph{octonionic K\"ahler} the parallel even Clifford structure defined by the vector bundles $E^{10},E^{12},E^{16},E^9$ on the Cayley-Rosenfeld projective planes $\EIII, \, \EVI, \, \EVIII , \, \FII$.

In conclusion, and with the exceptions of $$\EI = \frac{\Esix}{\Sp{4}}, \quad \EIV= \frac{\Esix}{\Ffour},\quad \EV=\frac{\Esix}{\SU{8}},$$ nine of the twelve exceptional Riemannian symmetric spaces of compact type admit at least one parallel even Clifford structure. 

Aim of the present paper is to describe how, basing on the recent work \cite{pp1}, \cite{pp3}, \cite{ppv} about $\Spin{9}$, $\Spin{10}\cdot\U{1}$ and further even Clifford structures, one can construct canonical differential 8-forms on the symmetric spaces $\EIII, \, \EVI, \, \EVIII , \, \FII$. Their classes are one of the cohomology generators, namely the one corresponding to their octonionic K\"ahler structure. I will discuss in particular for which of our nine exceptional Riemannian symmetric spaces of compact type the de Rham cohomology is fully canonical, i. e. fully generated by classes represented by canonical forms associated with parallel even Clifford structures.

\vspace{0.5cm}

\section{Poincar\'e polynomials}\label{poincare}

The following Table A collects some informations on the exceptional symmetric spaces of compact type. For each of them the dimension, the existence of torsion in the integral cohomology, the K\"ahler or quaternion K\"ahler or octonionic K\"ahler (K/qK/oK) property, the Euler characteristic $\chi$, and the Poincar\'e polynomial (up to mid dimension) are listed. The last column contains the references where the above informations are taken from. The Euler characteristics can be confirmed via the theory of elliptic genera, cf. \cite{hs}.

\begin{table}[H]
\caption{Exceptional symmetric spaces of compact type}
\begin{center}
\renewcommand{\arraystretch}{2.30}
\resizebox*{1.00\textwidth}{!}{
\begin{tabular}{|c||c|c|c|c|c|c|}
\hline
 & dim &torsion&K/qK/oK&$\chi$&Poincar\'e polynomial $P(t)=\sum_{i=0, \dots} b_i t^i$&Reference\\
\hline\hline
$\EI$&42&yes&&4&$1+t^8+t^9+t^{16}+t^{17}+t^{18}+ \dots $&\cite{is}\\
\hline
$\EII$&40&yes &qK &36& $1+t^4+t^6+2t^{8}+t^{10}+3t^{12}+ 2t^{14}+3t^{16}+2t^{18}+4t^{20}+\dots $&\cite{is2}\\
\hline
$\EIII$&32& no& K/oK &27&  $1+t^2+t^4+t^6+2(t^{8}+t^{10}+t^{12}+ t^{14})+3t^{16}+\dots $&\cite{tw}\\
\hline
$\EIV$&26& no&&0&$1+t^9+\dots$&\cite{ar}\\
\hline
$\EV$&70&yes &&72&$1+t^6+t^8+t^{10}+t^{12}+2(t^{14}+ t^{16}+t^{18}+ t^{20})+ \qquad \qquad $ &\\
&&&&&$\qquad \qquad +3(t^{22}+t^{24}+t^{26}+t^{28})+4(t^{30}+t^{32})+3t^{34}+\dots $ &\\
\hline
$\EVI$&64& yes& qK/oK&63&$1+t^4+2t^8+3t^{12}+4t^{16}+5t^{20}+6(t^{24}+t^{28})+7t^{32}+\dots $&\cite{n} \\
\hline
$\EVII$&54&no&K&56&$1+t^2+t^4+t^6+t^8+2(t^{10}+t^{12}+t^{14}+ t^{16})+ \qquad \qquad$&\cite{w}\\
&&&&&$\qquad \qquad +3(t^{18}+ t^{20}+t^{22}+t^{24}+t^{26})+\dots $&\\
\hline
$\EVIII$&128 &yes&oK&135&$1+t^8+t^{12}+2(t^{16}+ t^{20})+3(t^{24}+t^{28})+5t^{32}+\qquad \qquad$ &\\
&&&&&$ \qquad \qquad +4t^{36}+6(t^{40}+t^{44})+7(t^{48}+t^{52})+8t^{56}+7t^{60}+9t^{64}+\dots $&\\
\hline
$\EIX$&112& yes&qK& 120&$1+t^4+t^8+2(t^{12}+t^{16})+3 t^{20}+4(t^{24}+t^{28})+ \qquad \qquad $&\cite{s}\\
&&&&&$\qquad \qquad +5t^{32}+6(t^{36}+t^{40})+7(t^{44}+t^{48}+t^{52})+8t^{56}+\dots $&\\
\hline
$\FI$&28 & yes&qK&12&$1+t^4+2(t^8+t^{12})+\dots$&\cite{it} \\
\hline
$\FII$&16 &no&oK&3&$1+t^8+\dots$&\cite{bh} \\
\hline
$\GI$&8 & yes& qK&3&$1+t^4+\dots$ &\cite{bh}\\
\hline
\end{tabular}
}
\end{center}
\end{table}

Most of the de Rham cohomology structures are in the literature, according to the mentioned references. The author was not able to find a reference for the cohomology computations of $\EV$ and $\EVIII$. Their de Rham cohomologies can however be obtained through the following Borel presentation, cf. the original A. Borel article \cite{bo}, as well as \cite{n}.

Let $G$ be a compact connected Lie group, let $H$ be a closed connected subgroup of $G$ of maximal rank, and $T$ a common maximal torus. The de Rham cohomology of $G/H$ can be computed in terms of those of the classifying spaces $BG$, $BH$, $BT$ according to
\[
H^*(G/H) \stackrel{\sim}\longleftarrow H^*(BH)/\rho^* H^+(BG) \cong H^*(BT)^{W(H)}/(H^+(BT)^{W(G)}),
\]
thus as quotient of the ring $H^*(BT)^{W(H)}$ of invariants of the Weyl group $W(H)$. Here notations refer to the fibration $$G/H  \stackrel{i}\rightarrow BH \stackrel{\rho}\rightarrow BG.$$ 
Also $H^+= \oplus_{i>0} H^i$, and $(H^+(BT)^{W(G)})$ is the ideal of $H^*(BT)^{W(H)}$ generated by $H^+(BT)^{W(G)}$. 

The two cases of interest for us are 
\[
\EV: \; (G,H)= (\mathrm{E}_7, \SU{8}) \quad \text{and} \quad
\EVIII: \;  (G,H)= (\mathrm{E}_8, \mathrm{Spin}^+ (16)).
\]
In fact, the rings of invariants of the Weyl groups $W(\mathrm{E}_7)$, $W(\mathrm{E}_8)$ have been computed in \cite{tw}, \cite{n}, \cite{n2}. They read:
\[
H^*(BT)^{W(\mathrm{E}_7)} \cong \mathbb R[\sigma_2, \sigma_6, \sigma_{8}, \sigma_{10}, \sigma_{12}, \sigma_{14},\sigma_{18}], 
\]
\[
H^*(BT)^{W(\mathrm{E}_8)} \cong \mathbb R[\rho_2, \rho_{8}, \rho_{12},\rho_{14},\rho_{18},\rho_{20},\rho_{24},\rho_{30}],
\]
with $\sigma_\beta, \rho_\beta \in H^{2\beta}$. 

As mentioned, $\EV$ and $\EVIII$ are quotients respectively of $\mathrm{E}_7$  by $\SU{8}$ and of $\mathrm{E}_8$ by the subgroup $\mathrm{Spin^+}(16)$, a $\ZZ_2$ quotient of $\Spin{16}$ that is not $\mathrm{SO}(16)$, see \cite{it}. Since $\SU{8}$ has the same rank 7 of $\mathrm{E}_7$ and 
$\mathrm{Spin^+}(16)$ has the same rank 8 of $\mathrm{E}_8$, we can use Borel presentation. 

With Chern, Euler and Pontrjagin classes notations:
\[
H^*(BT)^{W(\mathrm{SU}(8))} \cong \mathbb R[c_2,c_3,c_4,c_5,c_6,c_7,c_8], 
\]
\[
H^*(BT)^{W(\mathrm{Spin^+}(16))} \cong \mathbb R[e,p_1,p_2,p_3,p_4,p_5,p_6,p_7],
\]
where $c_\alpha \in H^{2\alpha}, e \in H^{16}, p_\alpha \in H^{4\alpha}$. The cohomologies of $\EV$ and $\EVIII$ are now easily obtained. By interpreting the $\sigma_\beta$ and $\rho_{\beta}$ as relations among polynomials in the mentioned Chern, Euler and Pontrjagin classes, the Poincar\'e polynomials of $\EV$ and $\EVIII$ included in Table A follow from a straightforward computation.

Next, the following Table B contains the primitive Poincar\'e polynomials $\widetilde{P}(t)=\sum_{i=0, \dots} \widetilde{b}_i t^i$ of the nine exceptional Riemannian symmetric spaces that admit an even parallel Clifford structures. Here "primitive" has a different meaning, according to the considered  K/qK/oK structure.  Thus, for the Hermitian symmetric spaces $\EIII$ and $\EVII$, they are simply the polynomials with coefficients the primitive Betti numbers $$\widetilde{b}_i  = \dim \; (\ker [L_\omega^{n-i+1}: H^i \rightarrow H^{2n-i+2}]),$$ where $L_\omega$ is the Lefschetz operator, the multiplication of cohomology classes with that of complex K\"ahler form $\omega$, and $n$ is the complex dimension. 

In the positive quaternion K\"ahler setting, one has the vanishing of odd Betti numbers and the injectivity of the Lefschetz operator $L_\Omega: H^{2k-4} \rightarrow H^{2k}$, $k \leq n$, now with $\Omega$ the quaternion 4-form and $n$ the quaternionic dimension. A remarkable aspect of the primitive Betti numbers $$\widetilde{b}_{2k} = \dim ( \text{coker}  [L_\Omega: H^{2k-4} \rightarrow H^{2k}])$$  for positive quaternion K\"ahler manifolds is their coincidence with the ordinary Betti numbers of the associated Konishi bundle, the 3-Sasakian manifold fibering over it, cf. \cite{ga-sa}, p. 56. Indeed, one can check this coincidence on the exceptional Wolf spaces: just compare the quaternion K\"ahler part of Table B with Table III in the quoted  paper by Krzysztof Galicki and Simon Salamon. 

Finally, on the four Cayley-Rosenfeld planes, one still has the vanishing of odd Betti numbers and the injectivity of the map $$L_\Phi: H^{2k-8} \rightarrow H^{2k},$$ defined by multiplication with the octonionic 8-form $\Phi$, and with $k \leq 2n$, $n$ now the octonionic dimension, cf. Sections \ref{8-form} and \ref{CR}.

\begin{table}[H]
\caption{Primitive Poincar\'e polynomials $\widetilde{P}(t)=\sum_{i=0, \dots} \widetilde{b}_i t^i$}
\begin{center}
\renewcommand{\arraystretch}{2.30}
\resizebox*{1.0\textwidth}{!}{
\begin{tabular}{|c||c|}
\hline
Hermitian symmetric spaces&K\"ahler primitive Poincar\'e polynomial\\
\hline\hline
$\EIII$&  $1+t^{8}+t^{16}$\\
\hline
$\EVII$&$1+t^{10}+t^{18}$\\
\hline\hline\hline
 Wolf spaces& Quaternion K\"ahler primitive Poincar\'e polynomial \\
\hline\hline
$\EII$& $1+t^6+t^{8}+t^{12}+ t^{14}+t^{20} $\\
\hline
$\EVI$&$1+t^8+t^{12}+t^{16}+t^{20}+t^{24}+t^{32}$\\
\hline
$\EIX$&$1+t^{12}+t^{20}+t^{24}+t^{32}+t^{36}+t^{44}+t^{56} $\\
\hline
$\FI$&$1+t^8$\\
\hline
$\GI$&$1$ \\
\hline\hline\hline
 Cayley-Rosenfeld planes&Octonionic K\"ahler primitive Poincar\'e polynomial \\
\hline\hline
$\EIII$& $1+t^2+t^4+t^6+t^{8}+t^{10}+t^{12}+ t^{14}+t^{16}$\\
\hline
$\EVI$&$1+t^4+t^8+2(t^{12}+t^{16}+t^{20})+3(t^{24}+t^{28}+t^{32})$ \\
\hline
$\EVIII$&$1+t^{12}+t^{16}+ t^{20}+t^{24}+t^{28}+t^{32}+t^{36}+t^{40}+t^{44}+t^{48}+t^{52}+t^{56}+t^{60}+t^{64}$\\
\hline
$\FII$&$1$ \\
\hline
\end{tabular}
}
\end{center}
\end{table}

Note that $\EIII$ and $\EVI$ appear twice in the former table. The intersection of their primitive cohomology with respect with to the K/ok, respectively qK/oK structure, give rise to the following "fully primitive Poincar\'e polynomials", listed in the following table for the nine exceptional symmetric spaces of compact type admitting an even Clifford structure.

\begin{table}[H]
\caption{Fully primitive Poincar\'e polynomials $\overline{P}(t)=\sum_{i=0, \dots} \overline{b}_i t^i$}
\begin{center}
\renewcommand{\arraystretch}{2.30}
\resizebox*{1.0\textwidth}{!}{
\begin{tabular}{|c||c|}
\hline
%
%
Even Clifford symmetric spaces &Fully primitive Poincar\'e polynomial \\
\hline\hline
$\EII$& $1+t^6+t^{8}+t^{12}+ t^{14}+t^{20} $\\
\hline
$\EIII$&  $1 $\\
\hline
$\EVI$&$1+t^{12}+t^{24}$\\
\hline
$\EVII$&$1+t^{10}+t^{18}$\\
\hline
$\EVIII$&$1+t^{12}+t^{16}+ t^{20}+t^{24}+t^{28}+t^{32}+t^{36}+t^{40}+t^{44}+t^{48}+t^{52}+t^{56}+t^{60}+t^{64}$\\
\hline
$\EIX$&$1+t^{12}+t^{20}+t^{24}+t^{32}+t^{36}+t^{44}+t^{56} $\\
\hline
$\FI$&$1+t^8$\\
\hline
$\FII$&$1$ \\
\hline
$\GI$&$1$ \\
\hline
\end{tabular}
}
\end{center}
\end{table}

\vspace{0.5cm}

\section{The octonionic K\"ahler 8-form}\label{8-form}

An even Clifford structure on the Cayley-Rosenfeld projective planes $\FII ,$ $ \EIII , $ $ \EVI  , $ $  \EVIII$ is given by a suitable vector sub-bundle of their endomorphism bundle. 

To describe these vector sub-bundles look first, as proposed by Th. Friedrich \cite{fr}, at the following matrices, defining self-dual anti-commuting involutions in $\mathbb R^{16}$. It is natural to name them the \emph{octonionic Pauli matrices}:

\begin{equation*}\label{Pauli}
\begin{split}
S_0=\left(
\begin{array}{r|r}
0 & \Id \\
\hline 
\Id & 0
\end{array}\right), \quad
&S_1=\left(
\begin{array}{c|c}
0 & -R_i \\
\hline 
R_i& 0
\end{array}\right), \quad 
S_2=\left(
\begin{array}{c|c}
0 & -R_j \\
\hline 
R_j& 0
\end{array}\right), \\
S_3=\left(
\begin{array}{c|c}
0 & -R_k\\
\hline 
R_k& 0
\end{array}\right), \quad 
&S_4=\left(
\begin{array}{c|c}
0 & -R_e \\
\hline 
R_e& 0
\end{array}\right), \quad
S_5=\left(
\begin{array}{r|r}
0 & -R_f \\
\hline 
R_f & 0
\end{array}\right), \\ 
S_6=\left(
\begin{array}{c|c}
0 & -R_g \\
\hline 
R_g& 0
\end{array}\right), \quad 
&S_7=\left(
\begin{array}{c|c}
0 & -R_h \\
\hline 
R_h& 0
\end{array}\right), \quad
S_8=\left(
\begin{array}{r|r}
\Id & 0 \\
\hline
0 & -\Id
\end{array}\right).
\end{split}
\end{equation*}

\noindent Here $R_u$ is the right multiplication by the unit basic octonion $u=i,j,k,e,f,g,h$, and matrices $S_\alpha$ act on $\OO^2 \cong \RR^{16}$. As mentioned, for all $\alpha =0, \dots , 8$ and $\alpha \neq \beta$, one has
\[
S^*_\alpha = S_\alpha, \qquad  S^2_\alpha = \Id, \qquad S_\alpha S_\beta = - S_\beta S_\alpha.
\]

Next, on the real vector spaces $$\CC^{16}= \CC \otimes \RR^{16}, \qquad \HH^{16} =\HH \otimes \RR^{16}, \qquad \OO^{16}=\OO \otimes \RR^{16},$$ besides the endomorphisms $S_\alpha$ ($\alpha = 0, \dots , 8)$, thought now acting on the factor $\RR^{16}$, look also at the skew-symmetric endomorphisms:
\[
\mathfrak I \; \; \text{on} \; \; \CC^{16}, \qquad \mathfrak I , \, \mathfrak J , \, \mathfrak K \; \; \text{on} \; \; \HH^{16}, \qquad \mathfrak I , \, \mathfrak J , \,  \mathfrak K ,\,  \mathfrak E,\, \mathfrak F, \, \mathfrak G , \, \mathfrak H \; \; \text{on} \; \; \OO^{16},
\]
the multiplication on the left factor of the tensor product by the basic units in $\CC$, $\HH$, $\OO$. By enriching the nine $S_\alpha$ with such complex structures, we generate real vector subspaces $$E^{10} \subset \; \text{End} \; (\CC^{16}), \qquad E^{12} \subset \; \text{End} \; (\HH^{16}), \qquad E^{16} \subset\; \text{End} \; (\OO^{16}).$$ 

In order to get a Clifford map $$\varphi: \; \text{Cl}^0(E) \rightarrow \End{- \otimes \RR^{16}},$$ assume all generators to be anti-commuting with respect to a product formally defined as $\mathfrak I \wedge S_\alpha= - S_\alpha \wedge \mathfrak I, \dots$. Of course, the products $S_\alpha \wedge S_\beta$ and $\mathfrak I \wedge \mathfrak J, \dots$ are the usual compositions of endomorphisms.

It is convenient to use the notations:
\[
S_{-1} = \mathfrak I, \, S_{-2} = \mathfrak J, \, S_{-3} = \mathfrak K, \, S_{-4} = \mathfrak E, \, S_{-5} = \mathfrak F, \, S_{-6} = \mathfrak G, \, S_{-7} = \mathfrak H, 
\]
allowing to exhibit the mentioned even Clifford structures as
\[
E^9=<S_0, \dots , S_8>, \quad E^{10} =<S_{-1}> \oplus <S_0, \dots , S_8>,  
\]
\[ 
E^{12} =<S_{-1},S_{-2},S_{-3}> \oplus <S_0, \dots , S_8>,  
\]
\[
E^{16} = <S_{-1},\dots,S_{-7}> \oplus <S_0, \dots , S_8>,
\]
respectively defined on $\RR^{16}, \CC^{16}, \HH^{16}, \OO^{16}$. Note that the first summands of $E^{10},$ $E^{12},$ $E^{16}$ are the complex, the quaternionic and the octonionic structure in these linear spaces.

On Riemannian manifolds $M$, like the symmetric spaces $\FII ,$ $ \EIII , $ $ \EVI  , $ $  \EVIII$, the even Clifford structures are defined as vector bundles 
\[
E^9, \qquad E^{10} =E^1\oplus E^9, \qquad E^{12}=E^3 \oplus E^9, \qquad E^{16}=E^7 \oplus E^9,
\] 
with the line bundle $E^1=<S_{-1}>$ trivial, and $E^3$, $E^7$, $E^9$ locally generated by the complex structures $S_\alpha$ with negative index $\alpha$, and with the mentioned properties.

To allow a uniform notation, use the lower bound index $$A=0,-1,-3,-7,$$ according to whether $M=\FII , \EIII , \EVI  ,  \EVIII$, so that the mentioned generators of the even Clifford structure can be written as $$\{S_\alpha\}_{A \leq S_\alpha\leq 8}.$$ In all four cases one has a matrix of local almost complex structures
\[
J= \{J_{\alpha \beta}\}_{A \leq \alpha, \beta \leq 8},
\]
where $J_{\alpha \beta}= S_\alpha \wedge S_\beta$, so that $J$ is skew-symmetric. It is easily recognized that on the model linear spaces $\RR^{16}, \CC^{16}, \HH^{16}, \OO^{16}$, the upper diagonal elements $\{J_{\alpha \beta}\}_{A \leq \alpha < \beta \leq 8}$ are a basis of the Lie algebras $$\liespin{9} \subset \lieso{16}, \quad \liespin{10} \subset \liesu{16}, \quad \liespin{12} \subset \liesp{16}, \quad\liespin{16}.$$

Next, look at the skew-symmetric matrix of the (local) associated K\"ahler 2-forms:
\[
\psi= \{\psi_{\alpha \beta}\}_{A \leq \alpha, \beta \leq 8},
\]
and note that by the invariance property, the coefficients of its characteristic polynomial
\[
\det(tI-\psi)=t^{A+9}+\tau_2(\psi)t^{A+7}+\tau_4(\psi)t^{A+5} + \dots 
\]
give rise to global differential forms on $\FII ,$ $ \EIII , $ $ \EVI  , $ $  \EVIII$.

\begin{de} For $A=0,-1,-3,-7$, that is on the linear spaces $\RR^{16}, \CC^{16}, \HH^{16},\OO^{16}$, and on symmetric spaces  $\FII ,$ $ \EIII , $ $ \EVI  , $ $  \EVIII$, we call $$\Phi=\tau_4(\psi)$$ the \emph{octonionic K\"ahler 8-form}.
\end{de}
\vspace{0.5cm}

\section{Cayley-Rosenfeld planes}\label{CR}

The following statement has been proved in \cite{pp1}, \cite{pp3}; see also \cite{pp2}, \cite{oppv}.

\begin{te}\label{teo1}
(a) On the Cayley projective plane $\FII$, assuming the lower index $A=0$ in the matrix $\psi$, one has $\tau_2(\psi)=0$. Moreover, the octonionic K\"ahler 8-form
\[
\spinform{9} =\tau_4(\psi)
\]
is closed and its cohomology class generates the cohomology ring $H^*(\FII)$. 

(b) On the Hermitian symmetric space $\EIII$ one has (assuming now the lower index $A=-1$) $\tau_2(\psi)=-3\omega^2$, where $\omega$ is the complex K\"ahler 2-form of $\EIII$. The octonionic K\"ahler 8-form
\[
\spinform{10} =\tau_4(\psi)
\]
is closed and its class generates the K\"ahler primitive cohomology ring of $\EIII$. 

\end{te}

Another algebraic approach to the 8-form $\spinform{9}$, equivalent to the one outlined in the previous Section, has been proposed by  M. Castrill\'on L\`opez, P. M. Gadea and I. V. Mykytyuk, cf. \cite{cgm1}, \cite{cgm2}.

\begin{re}\label{quaternionic}
The above construction of the octonionic K\"ahler 8-form $\spinform{9}$ can be seen as parallel with the following way of looking at the quaternion K\"ahler 4-form $\Omega$ on quaternion Hermitian manifolds $M^8$. For this, consider the following \emph{quaternionic Pauli matrices}:

\begin{equation*}
\begin{split}
S^{\HH}_0=\left(
\begin{array}{r|r}
0 & \Id \\
\hline 
\Id & 0
\end{array}\right), \quad 
S^{\HH}_1=\left(
\begin{array}{c|c}
0 & -R^{\HH}_i \\
\hline 
R^{\HH}_i& 0
\end{array}\right), & \quad 
S^{\HH}_2=\left(
\begin{array}{c|c}
0 & -R^{\HH}_j \\
\hline 
R^{\HH}_j& 0
\end{array}\right), \\
S^{\HH}_3=\left(
\begin{array}{c|c}
0 & -R^{\HH}_k\\
\hline 
R^{\HH}_k& 0
\end{array}\right), \quad 
S^{\HH}_4 & =\left(
\begin{array}{r|r}
\Id & 0 \\
\hline
0 & -\Id
\end{array}\right),
\end{split}
\end{equation*}
acting on $\HH^2 \cong \RR^8$, where $R^{\HH}_u$ is the right multiplication by the unit basic quaternion $u=i,j,k$. As before, for $\alpha =0, \dots , 4$ and $\alpha \neq \beta$:
\[
(S^{\HH}_\alpha )^*= S^{\HH}_\alpha, \qquad  (S^{\HH}_\alpha )^2= \Id, \qquad S^{\HH}_\alpha S^{\HH}_\beta = - S^{\HH}_\beta S^{\HH}_\alpha.
\]
This can be applied to some Riemannian manifolds $M^8$, and the symmetric space $\GI = \mathrm G_2/\SO{4}$ is an example, by defining over it a Euclidean vector bundle $E^5$ locally spanned by five involutions with such properties, thus giving an even Clifford structure. This is equivalent to a quaternion Hermitian structure on $M^8$, as one can see by looking at the skew-symmetric matrix
\[
J= \{J_{\alpha \beta}\}_{0 \leq \alpha, \beta \leq 4}
\]
of almost complex structures $J_{\alpha \beta}= S^{\HH}_\alpha \circ S^{\HH}_\beta$, and at its associated matrix $\varphi=(\varphi_{\alpha\beta})$ of K\"ahler 2-forms. In the characteristic polynomial
\[
\det(tI-\varphi)=t^{5}+\tau_2(\varphi)t^{3}+\tau_4(\varphi)t,
\]
the coefficient $\tau_2(\varphi) \in \Lambda^4$ is easily seen to coincide with
\[
-2\Omega_L =-2[\omega^2_{\LH_i}+\omega^2_{\LH_j}+\omega^2_{\LH_k}],
\]
where $\Omega_L $ is the left quaternionic $4$-form, cf. \cite{pp1}, p. 329.

\end{re}

Going back to Theorem \ref{teo1}, some words of comment on the analogy between the definitions of $$\Phi_{\Spin{9}} \in \Lambda^8 (\FII) \; \quad \text{and} \; \quad \Phi_{\Spin{10}} \in \Lambda^8(\EIII).$$ When applied to the $\Spin{9} \subset \SO{16}$ and $\Spin{10} \subset \SU{16}$ structures on the linear spaces $\RR^{16}$ and $\CC^{16}$, both 8-forms can be written in terms of the cartesian coordinates. For example, the computation of $\Phi_{\Spin{9}}$ in $\RR^{16}$ gives a sum of 702 non zero monomials in the $dx_\alpha$ $(\alpha=1, \dots , 16)$ with coefficients $\pm 1, \pm 2, - 14$ (cf. \cite{pp1}, pp. 339-343 for the full description). 

One sees in particular that $\Phi_{\Spin{9}} \in \Lambda^8 \RR^{16}$ is, up to a constant, the 8-form defined by integrating the volume of octonionic lines in the octonionic plane. Namely, if $\nu_l$ denotes the volume form on the line $l = \{(x,mx)\}$ or $l = \{(0,y)\}$ in $\OO^2$,  a computation shows that
\[
\Phi_{\Spin{9}} (\RR^{16}) = c \int_{\OO P^1} p_l^*\nu_l \,dl,
\]
where $p_l: \OO^2 \cong \RR^{16} \rightarrow l$ is the orthogonal projection, $\OO P^1 \cong S^8$ is the octonionic projective line of all the lines $l \subset \OO^2$ and the constant $c$ turns out to be $\frac{39916800}{\pi^4}$, cf. \cite{pp1}, pp. 338-339. The integral in the right hand side of the previous formula is the definition of the octonionic 8-form in $\OO^2$ proposed by M. Berger \cite{BerCCP}. Of course, the octonionic lines of $\OO^2$ are distinguished 8-planes in $\RR^{16}$, so that this 1972 definition anticipates the spirit of calibrations in this context.

Now, in spite of the analogy between our constructions of $\Phi_{\Spin{9}}$ and $\Phi_{\Spin{10}}$, it is clear that a similar approach is not possible for $\Phi_{\Spin{10}} \in \Lambda^8(\CC^{16})$, due to the lack of a Hopf fibration to refer to. 

The following homological interpretation of $\Phi_{\Spin{10}}$ on $\EIII$ relates to its projective algebraic geometry, cf. \cite{pp3}.

\begin{pr}\label{main} Look at $\EIII$ as the closed orbit of the action of $\mathrm{E}_6$ in the projectified Jordan algebra of $3 \times 3$ Hermitian matrices over complex octonions, i. e. as the fourth Severi variety in $\CC P^{26}$ mentioned in the Introduction. Then the de Rham dual of the basis represented in $H^8(\EIII; \ZZ)$ by the forms $(\frac{1}{(2\pi)^4} \Phi_{\Spin{10}}, \frac{1}{(2\pi)^4} \omega^4)$ is given by the pair of algebraic cycles
\[
\Big(\CC P^4 + 3(\CC P^4)', \; \;  \CC P^4 + 5(\CC P^4)' \Big),
\]
where $\CC P^4, (\CC P^4)'$ are maximal linear subspaces, belonging to the two different families ruling a totally geodesic non-singular quadric $Q_8$ contained in $\EIII \subset \CC P^{26}$ as a complex octonionic projective line. 
\end{pr}

The described construction of the 8-form $\Phi$ has similarly properties on the remaining Cayley-Rosenfeld planes.

\begin{te}\label{teo2} 
(a) On the quaternion K\"ahler symmetric space $\EVI$, assuming the lower index $A=-3$ for the matrix $\psi$, the second coefficient $\tau_2(\psi)$ is proportional to the quaternion K\"ahler 4-form $\Omega$. Also, the octonionic K\"ahler 8-form
\[
\spinform{12} =\tau_4(\psi)
\]
is closed and its class is one of the two generators of the quaternion K\"ahler primitive cohomology ring of $\EVI$. 

(b) On the symmetric space $\EVIII$ one has again, assuming now $A=-7$, $\tau_2(\psi)=0$. The octonionic K\"ahler 8-form
\[
\spinform{16} =\tau_4(\psi)
\]
is closed and its class is one of the four generators of cohomology ring of $\EVIII$.

\end{te}

The closeness of the octonionic K\"ahler 8-form $\Phi = \tau_4(\psi)$ on the symmetric spaces $\FII,$ $ \EIII,$ $\EVI,$ $\EVIII$ is recognized by looking at the $\psi_{\alpha \beta}$ $(A \leq \alpha < \beta \leq 9$ and $A=0,-1,-3,-7$), as local curvature forms of a metric connection on the vector bundle defining their non-flat even Clifford structure. This allows, in the proof of Proposition \ref{main}, to relate the class of  $\tau_4(\psi)$ with the second Pontrjagin class of this bundle, see \cite{pp3} for details.

\begin{re} Note in Tables B and C the different behavior of the primitive cohomology in the four Cayley-Rosenfeld projective planes. Concerning $\FII, \EIII$, there is no "fully primitive cohomology" (and no torsion): the classes of $\spinform{9}$ in the first case and the classes of the K\"ahler 2-form $\omega$ and of the 8-form $\spinform{10}$ generate the whole cohomology. On the quaternion K\"ahler Wolf space there is (besides torsion) a unique primitive generator in $H^{12}$, whose description does not seem to follow our techniques. Finally $\EVIII$ has a much richer primitive cohomology with generators in $H^{12}, H^{16}, H^{20}$.
\end{re}
\vspace{0.5cm}

\section{Essential Clifford structures}

The discussion in the previous two Sections shows that the even Clifford structures allowing to get the octonionic K\"ahler form $\Phi = \tau_4(\psi)$ on $\EIII, \EVI,\EVIII$ follow a slightly different construction from that of $\FII$. The latter is in fact obtained by what is called a \emph{Clifford system} $C_m$, i. e. a vector sub-bundle of the endomorphisms bundle locally generated by $(m+1)$ self-dual anti-commuting involutions $S_\alpha$. For $\FII$ we chose $m+1=9$ and in the linear setting the $S_\alpha$ are the octonionic Pauli matrices. The parallel definition of the quaternionic Pauli matrices, outlined in Remark \ref{quaternionic}, gives a similar approach to $\Sp{2}\cdot\Sp{1}$ structures in dimension 8, using a Clifford system with $m+1=5$.

The vector bundles 

$$E^{10} \subset \, \text{End} \, (T \, \EIII), \quad E^{12} \subset \, \text{End} \, (T\, \EVI), \quad E^{16} \subset\, \text{End} \, (T\, \EVIII),$$ 
defining the even Clifford structure on the remaining Cayley-Rosenfeld planes, have both symmetric and skew-symmetric endomorphisms as local generators. This suggests the following:

\begin{de} We say that an even Clifford structure is \emph{essential} if it is not given as a Clifford system, i. e. if it cannot be locally generated by self-dual anti-commuting involutions. 
\end{de}

We have, cf. \cite{ppv}:

\begin{te} 

(i) Any even Clifford structure in dimension 64 and 128 is necessarily essential. In particular such are the parallel even Clifford structures on $\EVI$ and $\EVIII$.

(ii) The even Clifford structure on $\EIII$ is also essential, although in dimension 32 one can have a non-essential even Clifford structure.

\end{te}

The statement \emph{(i)} follows by a dimensional comparison of the representation spaces of pairs of Clifford algebras like ${\mathcal Cl}^0_{12}$ and ${\mathcal Cl}_{0,12}$, or ${\mathcal Cl}^0_{16}$ and ${\mathcal Cl}_{0,16}$. 

As for statement \emph{(ii)}, the same comparison of dimensions for the representation spaces of ${\mathcal Cl}^0_{10}$ and ${\mathcal Cl}_{0,10}$ leaves open the possibility of having a Clifford system with $m+1=10$ in the linear space $\RR^{32}$. In fact it is not difficult to write down such a Clifford system. However, the structure group of a $32$-dimensional manifold carrying such a Clifford system reduces to $\Spin{10} \subset \SU{16}$, and this would be the case of the holonomy group, assuming that such a Clifford system induces the parallel even Clifford structure of $\EIII$. Thus, $\EIII$ would have a trivial canonical bundle, in contradiction with the  positive Ricci curvature property of Hermitian symmetric spaces of compact type.
\vspace{0.5cm}

\end{document}